\documentclass[12pt]{article}
\usepackage{amsmath}\usepackage{amssymb}\usepackage{amscd}
\newtheorem{defi}{Definition}\newtheorem{theo}[defi]{Theorem}\newtheorem{prop}[defi]{Proposition}\newtheorem{lemm}[defi]{Lemma}

\def\R{\operatorname{R}}
\def\RR{\hat{\R}}
\def\h{{\mathbf{h}}}
\def\Diff{\text{Diff}_{\h}}
\def\der{\bar{\partial}} 
\def\Z{x} 
\def\th{\tilde{h}}
\def\q{\check{\operatorname{q}}}

\def\U{\operatorname{U}}
\def\Uh{\bar{\U}(\h)}
\def\gln{{\mathbf {gl}}_n}
\def\DRn{\mathcal{D}\!\,(\gln)}
\def\teL{\operatorname{L}}
\def\be{\begin{equation}}
\def\ee{\end{equation}}

\begin{document}
\begin{center}
{\Large\bf Rings of $\h$-deformed differential operators}

\vspace{.6cm}
{\bf B. Herlemont$^{\circ}$ and
O. Ogievetsky$^{\circ\ast}$\footnote{On leave of absence from P. N. Lebedev Physical Institute, Leninsky Pr. 53,
117924 Moscow, Russia} }

\vskip .4cm
$^\circ\ ${Aix Marseille Univ, Universit\'e de Toulon, CNRS, CPT, Marseille, France}\\
$^\ast${Kazan Federal University, Kremlevskaya 17, Kazan 420008, Russia}
\end{center}

\vskip .4cm
\hfill{\it In memory of Petr Kulish}

\vskip .4cm
\begin{abstract}

\vskip .2cm
\noindent We describe the center of the ring $\Diff(n)$ of $\h$-deformed differential operators of type A. 
We establish an isomorphism between certain localizations of $\Diff(n)$ 
and the Weyl algebra $\text{W}_n$ extended by $n$ indeterminates. 
\end{abstract}

\vskip .6cm
\section{Introduction} The ring $\Diff (n)$ of $\h$-deformed differential operators of type A appears in the theory of reduction algebras. 
A reduction algebra provides a tool to study decompositions of representations of an associative algebra with respect to its 
subalgebra in the situation when this subalgebra is the universal enveloping algebra of a reductive Lie algebra \cite{M,AST}.
We refer to \cite{T,Zh} for the general theory and uses of reduction algebras. 

\vskip .2cm 
Decompositions of tensor products of representations of a reductive Lie algebra $\mathfrak{g}$ is a particular case of a restriction problem, 
associated to the diagonal embedding of $\U(\mathfrak{g})$ into $\U(\mathfrak{g})\otimes\U(\mathfrak{g})$. The corresponding reduction algebra, 
denoted $\mathcal{D}\!\,(\mathfrak{g})$, is called ``diagonal reduction algebra" \cite{KO2}. 
A description of the diagonal reduction algebra $\DRn$ in terms of generators and (ordering) defining relations was given in \cite{KO2, KO3}. 

\vskip .2cm 
The diagonal reduction algebra $\DRn$ admits an analogue of the ``oscillator realization", in the rings $\Diff(n,N)$, $N=1,2,3,\dots$, of $\h$-deformed 
differential operators, see \cite{KO5}. The ring $\Diff(n, N)$ can be obtained by the reduction of the ring of differential
operators in $nN$ variables (that is, of the Weyl algebra $\text{W}_{nN}=\text{W}_{n}^{\otimes N})$ with respect to the natural action of $\gln$. Similarly to the ring of $q$-differential 
operators \cite{WZ}, the algebra $\Diff(n, N)$ can be described in the R-matrix formalism. The R-matrix, needed here, is a solution of the so-called 
dynamical Yang--Baxter equation (we refer to \cite{F,GN,ES} for different aspects of the dynamical Yang--Baxter equation and its solutions).   

\vskip .2cm 
The ring $\Diff(n,N)$ is formed by $N$ copies of the ring $\Diff(n)=\Diff(n,1)$. The aim of the present article is to investigate the structure of the 
ring $\Diff(n)$. Our first result is the description of the center of $\Diff(n)$: it is a ring of polynomials in $n$ generators. 

\vskip .2cm 
As follows from the results of \cite{KO4}, the ring $\Diff(n)$ is a noetherian Ore domain. It is therefore natural to investigate its field of fractions
and test the validity of the Gelfand--Kirillov-like conjecture \cite{GK}. For the ring of $q$-differential operators, the isomorphism (up to a certain localization and completion) with the Weyl algebra was given in \cite{O}. The second result of the present article consists in a construction of an isomorphism 
between certain localizations of $\Diff (n)$ and the Weyl algebra $\text{W}_n$ extended by $n$ indeterminates. In particular, our isomorphism implies the isomorphism of 
the corresponding fields of fractions.

\vskip .2cm 
According to the general theory of reduction algebras, the ring $\Diff(n)$ admits the action of Zhelobenko operators \cite{Zh,KO1} by automorphisms. 
The Zhelobenko operators generate the action of the braid group $\mathrm{B}_n$. However the Weyl algebra $\text{W}_n$ admits the action of the symmetric 
group $\mathrm{S}_n$ by automorphisms.
As a by-product of our construction we define the action of the symmetric group by automorphisms on the ring $\Diff(n)$. 
Moreover these formulas can be generalized to produce the action of the symmetric group by automorphisms on the rings $\Diff(n,N)$ for any $N$, 
and on the diagonal reduction algebra $\DRn$. 

\vskip .2cm 
Section \ref{hdiff-sect} contains the definition of the ring of $\h$-deformed differential operators and some of their properties used in the sequel. 
In Section \ref{quadce} we present a family of $n$ quadratic central elements. We then describe an $n$-parametric family of ``highest weight" 
representations of $\Diff(n)$ and calculate values of the quadratic central elements in these representations. In Section \ref{isowe} we introduce the necessary localizations
of $\Diff(n)$ and of the Weyl algebra, check 
the Ore conditions and establish the above mentioned isomorphism of the localized rings. In Section \ref{secco} we prove the completeness of the family 
of central elements constructed in Section \ref{quadce}. Then we describe the action of the symmetric group on $\Diff(n,N)$ and $\DRn$ as well as the action of 
the braid group, generated by Zhelobenko operators, on 
a localization of the Weyl algebra. Also, we present 
a $2n$-parametric family of representations of the algebra $\Diff(n)$ implied by our construction. 

\vskip .6cm\noindent
{\bf \large Notation}

\vskip .4cm\noindent
Throughout the paper, $\mathfrak{k}$ denotes the ground ring of characteristic zero. 

\vskip .2cm
The symbol $s_i$ stands for the transposition $(i,i+1)$. 

\vskip .2cm
We denote by $\U(\h)$ the free commutative $\mathfrak{k}$-algebra in generators $\th_i$, $i=1,\dots,n$. Set $\th_{ij}=\th_i-\th_j\in \h$. 
We define $\Uh$ to be the ring of fractions of 
$\U(\h)$ with respect to  the multiplicative set of denominators, generated by the elements
$(h_{ij}+k)^{-1}$, $k\in\mathbb{Z}$.
Let
\begin{equation}\label{defphichi}\psi_i:=\prod_{k:k>i}\th_{ik}\ , \psi_i':=\prod_{k:k<i}\th_{ik}\ \ \text{and}\ \ \chi_i:=\psi_i\psi_i'\ ,\ i=1,\dots,n\ .\end{equation}
Let $\varepsilon_j$, $j=1,\dots,n$, be the elementary translations of the generators of $\U(\h)$, 
$\varepsilon_j:\th_i\mapsto\th_i+\delta_i^j$. For an element $p\in \Uh$ we denote  
$\varepsilon_j(p)$ by $p[\varepsilon_j]$. 

\section{Definition and properties of rings of $\h$-deformed differential operators}\label{hdiff-sect}
Let $\RR=\{\RR_{ij}^{kl}\}_{i,j,k,l=1}^n$ be a matrix of elements of $\Uh$, with nonzero entries 
\begin{equation}\label{dynRcomp}
\RR_{ij}^{ij}=\frac{1}{\th_{ij}}\ ,\ \ i\not=j\ ,\qquad\text{and}\qquad \RR_{ji}^{ij}=\left\{
\begin{array}{cc}
 \dfrac{\th_{ij}^2-1}{\th_{ij}^2}\ ,& \, i<j,\\[0.5em]
 1\ ,&\, i\geq j\ .
\end{array}
\right. \end{equation}
The matrix $\RR$ is the standard solution of the dynamical Yang--Baxter equation
\begin{equation} \sum_{a,b,u}{\RR}^{ij}_{ab}\RR^{bk}_{ur}[-\varepsilon_a]\RR^{au}_{mn}=\sum_{a,b,u}
\RR^{jk}_{ab}[-\varepsilon_i]\RR^{ia}_{mu}
\RR^{ub}_{nr}[-\varepsilon_m]\ \end{equation}
of type A. 

\vskip .2cm
The ring $\Diff (n)$ of $\h$-deformed differential operators of type A is a $\Uh$-bimodule with the generators $\Z^j$ and $\der_j$, $j=1,\dots,n$.
The ring $\Diff (n)$ is free as a one-sided $\Uh$-module; the left and right $\Uh$-module structures are related by
\begin{equation}\label{relsweights}\th_i\Z^j=\Z^j(\th_i+\delta_i^j)\ ,\ \th_i\der_j=\der_j(\th_i-\delta_i^j)\ . \end{equation}
The defining relations for the generators $\Z^j$ and $\der_j$, $j=1,\dots,n$, read (see \cite{KO5})
\begin{equation}\label{defhdiffA}
 \Z^i \Z^j=\sum_{k,l}\RR_{kl}^{ij}\Z^k \Z^l\ ,\qquad \der_i\der_j=\sum_{k,l}\RR_{ji}^{lk}\der_k\der_l\ ,\qquad 
 \Z^i\der_j=\sum_{k,l}\RR_{lj}^{ki}[\varepsilon_k]\der_k \Z^l-\delta_{j}^i\ ,
\end{equation}
or, in components, 
\begin{equation}\label{relshdiffcompx}
\Z^i \Z^j=\frac{\th_{ij}+1}{\th_{ij}}\Z^j \Z^i\ ,\ i<j\ ,\end{equation}
\begin{equation}\label{relshdiffcompd}
\ \der_i \der_j=\frac{\th_{ij}-1}{\th_{ij}}\,\der_j \der_i\ ,\ i<j\ ,\end{equation}
\begin{equation}\label{relshdiffcompdx}
\Z^i \der_j=\left\{
\begin{array}{ll}
\der_j \Z^i&,\ \ i<j\ , \\[1em]
\displaystyle{\frac{\th_{ij}(\th_{ij}-2)}{(\th_{ij}-1)^2}}\,\der_j \Z^i&
,\  i>j\ ,\end{array}\right.\end{equation}
\begin{equation}\label{relshdiffcompdx2} 
\Z^i \der_i= \sum_j\displaystyle{\frac{1}{1-\th_{ij}}}\,\der_j \Z^j-1\ .
\end{equation}
The ring $\Diff (n)$ admits Zhelobenko automorphisms $\q_i$, $i=1,\dots,n-1$, given by (see \cite{KO5})
\begin{equation}\label{automq}\begin{split}
\q_i(\Z^i)&=-\Z^{i+1}\frac{\th_{i,i+1}}{\th_{i,i+1}-1},
\qquad \q_i(\Z^{i+1})=\Z^{i},\quad \q_i(\Z^j)=\Z^j,\quad j\not=i,i+1\ ,\\
\q_i(\der_{i})&=- \frac{\th_{i,i+1}-1}{\th_{i,i+1}}\der_{i+1},\qquad
\q_i(\der_{i+1})=\der_{i},\quad \q_i(\der_j)=\der_j,\quad
j\not=i,i+1\ ,\\[.5em]
\q_i(\th^j)&=\th_{s_i(j)}\ .
\end{split}
\end{equation}
The operators $\q_i$, $i=1,\dots,n-1$, generate the action of the braid group, see \cite{Zh,KO1}.

\vskip .2cm
The  ring $\Diff (n)$ admits an involutive anti-automorphism $\epsilon$, defined by 
\begin{equation}\label{antiautom}\epsilon(\th_i)=\th_i, \epsilon(\der_i)=\varphi_i\Z^i,\epsilon(\Z^i)=\der_i\varphi_i^{-1} ,
\, \text{where}\ \varphi_i:=\frac{\psi_i}{\psi_i[-\varepsilon_i]}=\prod_{k:k>i}\frac{\th_{ik}}{\th_{ik}-1},
\, i=1,\dots,n,\end{equation} 
The proof reduces to the formula 
$$\frac{\varphi_i[-\epsilon_j]}{\varphi_i}=\frac{\th_{ij}^2-1}{\th_{ij}^2}\ \ \text{for}\ \ 1\leq i<j\leq n\ .$$ 

The construction of central elements in the next Section uses the elements
\begin{equation}\label{defga} \Gamma_i:=\der_i\Z^i\ \ \text{for}\ \ i=1,\dots,n\ .\end{equation}
We collect some properties of these elements.
   
\begin{lemm}\hspace{-.2cm}. \label{qga}{\rm We have  
\begin{itemize}
\item[(i)] $\Gamma_i\Z^j=\displaystyle{\frac{\th_{ij}+1}{\th_{ij}}}\Z^j\Gamma_i\ $ and 
$\ \Gamma_i\der_j=\displaystyle{\frac{\th_{ij}-1}{\th_{ij}}}\der_j\Gamma_i\ $ for $i\neq j$, $\ i,j=1,\dots,n$.
\item[(ii)] $\q_i(\Gamma_j)=\Gamma_{s_i(j)}$ for $i=1,\dots,n-1$ and $j=1,\dots,n$.
\item[(iii)] $\Gamma_i\Gamma_j=\Gamma_j\Gamma_i\ $ for $i,j=1,\dots,n$.
\end{itemize}
}\end{lemm}
\noindent{\it{Proof.}} Formulas (i) and (ii) are obtained by a direct calculation; (iii) follows from (i). \hfill$\square$

\vskip .2cm
We will use the following technical Lemma. 
\begin{lemm}\hspace{-.2cm}.\label{lecommfamilies} {\rm 
Let $\mathfrak{A}$ be an associative algebra. Assume that elements $\breve{h}_i$, $\breve{Z}_i,\breve{Z}_i\in\mathfrak{A}$, $i=1,\dots,n$, satisfy
$$\breve{h}_i\breve{h}_j=\breve{h}_j\breve{h}_i\ ,\ \breve{h}_i\breve{Z}^j=\breve{Z}^j(\breve{h}_i+\delta_i^j)\ ,\ \breve{h}_i\breve{Z}_j=\breve{Z}_j(\breve{h}_i-\delta_i^j)\ ,\ 
i,j=1,\dots, n\ .$$
Let $\breve{h}_{ij}:=\breve{h}_i-\breve{h}_j$ and 
$$\breve{\psi}_i:=\prod_{k:k>i}\breve{h}_{ik}\ , \breve{\psi}_i':=\prod_{k:k<i}\breve{h}_{ik}\ ,\ i=1,\dots,n\ .$$
Assume that the elements $\breve{h}_{ij}$ are invertible. Then 

\vskip .2cm
\noindent (i) the elements $\breve{Z}^i$ satisfy $$\breve{Z}^i\breve{Z}^j=\frac{\breve{h}_{ij}+1}{\breve{h}_{ij}}\breve{Z}^j\breve{Z}^i\ \ \text{for}\ \  i<j\ ,\ i,j=1,\dots,n$$ if and only if 
any of the two families $\{\breve{Z}^{\circ i}\}_{i=1}^n$ or $\{\breve{Z}'^{\circ i}\}_{i=1}^n$ where 
\begin{equation}\label{twocommfamiliesa}  \breve{Z}^{\circ i}:=\psi_i\breve{Z}^i\ ,\ \breve{Z}'^{\circ i}:=\breve{Z}^i\psi_i'\end{equation}
is commutative; 

\vskip .2cm
\noindent (ii) the elements $\breve{Z}_i$ satisfy $$\breve{Z}_i\breve{Z}_j=\frac{\breve{h}_{ij}-1}{\breve{h}_{ij}}\breve{Z}_j\breve{Z}_i\ \ \text{for}\ \  i<j\ ,\ i,j=1,\dots,n$$ if and only if 
any of the two families $\{\breve{Z}^{\circ}_i\}_{i=1}^n$ or $\{\breve{Z}'^{\circ}_i\}_{i=1}^n$ where 
\begin{equation}\label{twocommfamiliesb}  \breve{Z}^{\circ}_i:=\psi_i\breve{Z}_i\ ,\ \breve{Z}'^{\circ}_i:=\breve{Z}_i\psi_i'\end{equation}
is commutative.}\end{lemm}
\noindent{\it{Proof.}} A direct calculation. \hfill$\square$

\section{Quadratic central elements}\label{quadce}
Let $e_k:=\sum_{i_1<\dots<i_k}\th_{i_1}\dots \th_{i_k}$, $k=0,\dots,n$, be the elementary symmetric functions in the variables $\th_1,\dots,\th_n$. Set
$$c_k:=\sum_j \frac{\partial e_k}{\partial\th_j}\Gamma_j-e_k\ ,$$
where $\Gamma_j$, $j=1,\dots,n$, are the elements defined in (\ref{defga}).

\vskip .2cm 
It follows from Lemma \ref{qga} that $\q_j(c_k)=c_k$ for all $j=1,\dots,n-1$ and $k=1,\dots,n$.

\begin{prop}\hspace{-.2cm}.\label{quceel} {\rm  The elements $c_k$, $k=1,\dots,n$, belong to the center of the ring $\Diff (n)$. }
\end{prop}
\noindent{\it{Proof.}} We shall use the generating functions
$$e(t):=\sum_{k=0}^ne_kt^k=\prod_i(1+\th_i t)$$ 
and 
$$c(t):=\sum_{k=1}^nc_kt^k=u(t)e(t)+1\ \ \text{with}\ \ u(t):=t\sum_i\frac{1}{1+\th_i t}\Gamma_i-1\ .$$
The expression $u(t)$ is introduced for convenience; the denominator $1+\th_i t$, which is not defined in the ring $\Diff(n)$, vanishes in the combination $u(t)e(t)$.

\vskip .2cm
We shall check that the polynomial $c(t)$ is central. We have 
\begin{equation}\label{quacen1}\Z^je(t)=\frac{1+(\th_j-1)t}{1+\th_jt}e(t)\Z^j\ .\end{equation}
Next, it follows from Lemma \ref{qga} that
$$\Z^ju(t)=\biggl(\sum_{k:k\neq j}\frac{t}{1+\th_k t}\frac{\th_{kj}}{\th_{kj}+1}
\Gamma_k+
\frac{t}{1+(\th_j-1)t}\Bigl(\sum_k\frac{1}{1-\th_{jk}}\Gamma_k-1\Bigr)-1\biggr)\Z^j\ .$$
The coefficient of $\Gamma_k$ in this expression is equal to $\frac{1+\th_j}{1+(\th_j-1)t}$ for both $k\neq j$ and $k=j$. Therefore, 
\begin{equation}\label{quacen2}\Z^ju(t)=\frac{1+\th_jt}{1+(\th_j-1)t}u(t)\Z^j\ .\end{equation}
Combining (\ref{quacen1}) and (\ref{quacen2}) we find that $c(t)$ commutes with $\Z^j$, $j=1,\dots,n$. For $\der_j$ one can either make a parallel calculation or use the anti-automorphism (\ref{antiautom}). \hfill$\square$

\begin{lemm}\hspace{-.2cm}. \label{gammaandc}{\rm  (i) The matrix $V$, defined by $V^k_j:=\frac{\partial e_j}{\partial \th_k}$, is invertible. 
Its inverse is 
$$(V^{-1})^j_i=\frac{(-1)^{j-1}\th_i^{n-j}}{\chi_i}\ ,$$
where the elements $\chi_i$ are defined in (\ref{defphichi}). 
 
\vskip .2cm\noindent 
(ii) We have
\begin{equation}\label{znachc}\chi_j\Gamma_j=\th_j^n-\th_j^nc(-\th_j^{-1})\ .\end{equation}
}
\end{lemm}

\vskip .2cm\noindent{\it{Proof.}} (i) See, e.g. \cite{OP}, Proposition 4.

\vskip .2cm\noindent
(ii) Rewrite the equality $c_k=\sum_jV_k^j\Gamma_j-e_k$ in the form
$$\Gamma_j=\sum_k(V^{-1})_j^k(c_k+e_k)=\frac{1}{\chi_j}\sum_k(-1)^{k-1}\th_j^{n-k}(c_k+e_k)=-\frac{\th_j^n}{\chi_j}(c(-\th_j^{-1})+e(-\th_j^{-1})-1)\ .$$
Since $e(-\th_j^{-1})=0$, we obtain (\ref{znachc}). \hfill$\square$

\paragraph{Highest weight representations.} The ring $\Diff(n)$ admits an $n$-parametric family of ``highest weight" representations. To define them, let $\mathfrak{D}_n$ be 
an $\Uh$-subring of $\Diff(n)$ generated by $\{\der_i\}_{i=1}^n$. 
Let $\vec{\lambda}:=\{\lambda_1,\dots,\lambda_n\}$ be a sequence of length $n$ of complex numbers such that $\lambda_i-\lambda_j
\notin\mathbb{Z}$ for all $i,j=1,\dots,n$, $i\neq j$. Denote by $M_{\vec{\lambda}}$ the 
one-dimensional $\mathfrak{k}$-vector space
with the basis vector $\vert\;\rangle$. Under the specified conditions on $\vec{\lambda}$ the formulas  
$$\th_i\colon \vert\; \rangle\mapsto\lambda_i\vert\;\rangle\ ,\ \der_i\colon \vert\;\rangle\mapsto 0\ ,\ i=1,\dots,n\ ,$$
define the $\mathfrak{D}_n$-module structure on $M_{\vec{\lambda}}$. 
We shall call the induced representation $\text{Ind}_{\mathfrak{D}_n}^{\Diff(n)}M_{\vec{\lambda}}$ the ``highest weight representation" of highest weight 
$\vec{\lambda}$. 

\begin{lemm}\hspace{-.2cm}. \label{valhwcel}{\rm The central operator $c_k$, $k=1,\dots,n$, acts on the module $\text{Ind}_{\mathfrak{D}_n}^{\Diff(n)}M_{\vec{\lambda}}$ by scalar multiplication 
on $-e_k\vert_{\th_i\mapsto \lambda_i-1}$, the evaluation of the symmetric function $-e_k$ on the shifted vector $\{\lambda_1-1,\dots,\lambda_n-1\}$.  }
\end{lemm}

\vskip .2cm\noindent{\it{Proof.}} It is sufficient to calculate the value of $c_k$ on the highest weight vector $\vert\;\rangle$. In terms of generating functions we 
have to check that 
$$e(t)u(t)\colon \vert\;\rangle\mapsto -\prod_i\left(1+(\lambda_i-1)t\right)\ \vert\;\rangle\ .$$
It follows from 
\cite{KO5}, section 3.3, that 
$$\Gamma_j\ \vert\;\rangle = \frac{\chi_j[\varepsilon_j]}{\chi_j}\ \vert\;\rangle\ ,\ j=1,\dots,n\ .$$
Therefore we have to check that 
\begin{equation}\label{valingefu}e(t)\left(t\sum_i\frac{1}{1+\th_i t}\frac{\chi_j[\varepsilon_j]}{\chi_j}-1\right)\colon \vert\;\rangle\mapsto -\prod_i\left(1+(\lambda_i-1)t\right)\ \vert\;\rangle\ .\end{equation}
We use another formula from \cite{KO5} (Note 3 after the proof of Proposition 4.3 in Section 4.2)
$$\prod_l\frac{\th_0-\th_l-1}{\th_0-\th_l}+\sum_j\frac{1}{\th_0-\th_j}\,\frac{\chi_j[-\varepsilon_j]}{\chi_j}=1\ ,$$
where $\th_0$ is an indeterminate. After the replacements $\th_0\to t^{-1}$ and $\th_j\to -\th_j$, $j=1,\dots,n$, this formula becomes 
$$\frac{e(t)[-\varepsilon]}{e(t)}+t\sum_j\frac{1}{1+\th_jt}\frac{\chi_j[\varepsilon_j]}{\chi_j}=1\ ,$$
where $\varepsilon=\varepsilon_1+\dots+\varepsilon_n$, which implies (\ref{valingefu}). \hfill$\square$

\section{Isomorphism between rings of fractions}\label{isowe}
It follows from the results of \cite{KO4} that the ring $\Diff (n)$ has no zero divisors. Let $\text{S}_\Z$ be the 
multiplicative set generated by $\Z^j$, $j=1,\dots,n$. The set $\text{S}_\Z$ satisfies both left and right Ore conditions (see, e.g., \cite{A} for definitions): say, for the 
left Ore conditions 
we have to check only that for any $\Z^k$ and a monomial $m=\der_{i_1}\dots \der_{i_A}\Z^{j_1}\dots \Z^{j_B}$ there exist $\tilde{s}\in\text{S}_\Z$ and $\tilde{m}\in\Diff(n)$ such that 
$\tilde{s}m=\tilde{m}\Z^k$. The structure of the commutation relations (\ref{relshdiffcompx}-\ref{relshdiffcompdx2}) shows that one can choose 
$\tilde{s}=(\Z^k)^{\nu}$ with sufficiently large $\nu$. Denote by $\text{S}_\Z^{-1}\Diff (n)$ the localization of the ring $\Diff (n)$ with respect to the set $\text{S}_\Z$. 

\vskip .2cm
Let $\text{W}_n$ be the Weyl algebra, the algebra with the generators $X^j,D_j$, $j=1,\dots,n$, and the defining relations
$$X^iX^j=X^jX^i\ ,\ D_iD_j=D_jD_i\ ,\ D_iX^j=\delta_i^j+X^jD_i\ ,\ \ i,j=1,\dots,n\ .$$
Let $\text{T}$ be the multiplicative set generated by $X^jD_j-X^kD_k+\ell$, $1\leq j<k\leq n$, $\ell\in\mathbb{Z}$, and $X^j$, $j=1,\dots,n$. The set $\text{T}$ 
satisfies left and right Ore conditions (see \cite{KO4}, Appendix).
Denote by $\text{T}^{-1}\text{W}_n$ the localization of $\text{W}_n$ relative to the set $\text{T}$.

\vskip .2cm
Let $a_1,\dots,a_n$ be a family of commuting variables. We shall use the following notation: 
$$\mathcal{H}_j:=D_jX^j\ ,\ \mathcal{H}_{jk}:=\mathcal{H}_j-\mathcal{H}_k\ ,$$ 
$$\Psi'_j:=\prod_{k:k<j}\mathcal{H}_{jk}\ ,\ \Psi_j:=\prod_{k:k>j}\mathcal{H}_{jk}\ ,$$
$$\mathbf{C}(t):=\sum_{k=1}^n a_kt^k\ ,\ \Upsilon_i:=\mathcal{H}_i^n \left(1-\mathbf{C}(-\mathcal{H}_i^{-1})\right)\ .$$  
The polynomial $\mathbf{C}$ has degree $n$ so the element $\Upsilon_i$ is a polynomial in $\mathcal{H}_i$, $i=1,\dots,n$.

\begin{theo}\hspace{-.2cm}. {\rm  The ring  $\text{S}_\Z^{-1}\Diff (n)$ is isomorphic to the ring $\mathfrak{k}[a_1,\dots,a_n]\otimes \text{T}^{-1}\text{W}_n$.}
\end{theo}
\noindent{\it{Proof.}} The knowledge of the central elements (Proposition \ref{quceel}) allows to exhibit a generating set of the ring $\text{S}_\Z^{-1}\Diff (n)$ in 
which the required isomorphism is quite transparent. 

\vskip .2cm
In the localized ring $\text{S}_\Z^{-1}\Diff (n)$ we can use the set of generators $\{\th_i,\Z^i,\Gamma_i\}_{i=1}^n$ instead of $\{\th_i,\Z^i,\der_i\}_{i=1}^n$. 
By Lemma \ref{gammaandc} (ii), $\{\th_i,\Z^i,c_i\}_{i=1}^n$ is also a generating set. Finally, $\mathfrak{B}_{\text{D}}:=\{\th_i,\Z^{\circ i},c_i\}_{i=1}^n$, 
where $\Z'^{\circ i}:=\Z^i\psi_i'\ ,\ i=1,\dots,n$, is a generating set of the localized ring 
$\text{S}_\Z^{-1}\Diff (n)$ as well. It follows from Lemma \ref{lecommfamilies} that the family $ \{\Z'^{\circ i}\}_{i=1}^n$ is commutative. The complete set of the 
defining relations for the generators from the set $\mathfrak{B}_{\text{D}}$ reads 
\begin{equation}\label{defrebd}\begin{array}{l}\th_i\th_j=\th_j\th_i\ ,\ \th_i \Z'^{\circ j}=\Z'^{\circ j}(\th_i+\delta_i^j)\ ,\ \Z'^{\circ i}\Z'^{\circ j}=\Z'^{\circ j}\Z'^{\circ i}\ ,\ i,j=1,\dots,n\ ,\\[.5em]
c_i\ \text{are central}\ ,\ i=1,\dots,n\ .\end{array}\end{equation}

In the localized ring $\mathfrak{k}[a_1,\dots,a_n]\otimes \text{T}^{-1}\text{W}_n$ we can pass to the set of generators $\mathfrak{B}_{\text{W}}:=\{\mathcal{H}_i,X^i,a_i\}_{i=1}^n$ with the defining relations
\begin{equation}\label{defrebw}\begin{array}{l}\mathcal{H}_i\mathcal{H}_j=\mathcal{H}_j\mathcal{H}_i\ ,\ \mathcal{H}_i X^j=X^j(\mathcal{H}_i+\delta_i^j)\ ,\ X^{i}X^{j}=X^j X^i\ ,\ i,j=1,\dots,n\ ,\\[.5em]
a_i\ \text{are central}\ ,\ i=1,\dots,n\ .\end{array}\end{equation}
The comparison of (\ref{defrebd}) and (\ref{defrebw}) shows that we have the isomorphism $$\mu\colon \mathfrak{k}[a_1,\dots,a_n]\otimes \text{T}^{-1}\text{W}_n\to 
\text{S}_\Z^{-1}\Diff (n)$$ 
given on our generating sets $\mathfrak{B}_{\text{D}}$ and $\mathfrak{B}_{\text{W}}$  by
\begin{equation}\label{fromwtodiffnb} \mu\colon  X^i\mapsto \Z'^{\circ i}\ ,\ \mathcal{H}_i\mapsto \th_i
\ ,\ a_i\mapsto c_i\ ,\ i=1,\dots,n\ .\end{equation}
The proof is completed.\hfill$\square$

\vskip .2cm
\noindent We shall now rewrite the formulas for the isomorphism $\mu$ in terms of the original generators of the rings $\text{S}_\Z^{-1}\Diff (n)$ and 
$\mathfrak{k}[a_1,\dots,a_n]\otimes \text{T}^{-1}\text{W}_n$. 

\begin{lemm}\hspace{-.2cm}. {\rm We have
\begin{equation}\label{fromwtodiff} \mu\colon  X^i\mapsto \Z^i\psi_i'\ ,\ D_i\mapsto (\psi_i')^{-1}\th_i(\Z^i)^{-1}\ ,\ a_i\mapsto c_i\ ,\ i=1,\dots,n\ .\end{equation}

\noindent and
\begin{equation}\label{fromdifftow} \mu^{-1}\colon  \th_i\mapsto \mathcal{H}_i\ ,\ \Z^i\mapsto X^i\frac{1}{\Psi_i'}\ ,\ \der_i\mapsto \frac{\Upsilon_i}{\Psi_i}
(X^i)^{-1}\ ,\ i=1,\dots,n\ .\end{equation}
}\end{lemm}

\noindent{\it{Proof.}} We shall comment only the last formula in (\ref{fromdifftow}). Lemma \ref{gammaandc} part (ii) implies that $\mu^{-1}(\chi_i\Gamma_i)=
\Upsilon_i$ and the formula for $\mu^{-1}(\der_i)$ follows since $\der_i=\Gamma_i(\Z^i)^{-1}$. \hfill$\square$

\section{Comments}\label{secco}
We shall now establish several corollaries of our construction. 

\paragraph{1.} We can now give the description of the center of the ring $\Diff (n)$.
\begin{lemm}\hspace{-.2cm}. {\rm The center of the ring 
$\Diff (n)$ is formed by polynomials in the elements $\{c_i\}_{i=1}^n$.}\end{lemm}

\noindent{\it{Proof.}}  This is a direct consequence of the defining relations (\ref{defrebd}) for the generating set $\mathfrak{B}_{\text{D}}$. 
\hfill$\square$
 
\paragraph{2.} The symmetric group $\mathrm{S}_n$ acts by automorphisms on the algebra $\text{W}_n$, 
$$\pi(X^j)=X^{\pi(j)}\ ,\ \pi(D_j)=D_{\pi(j)}\ \ \text{for}\ \  \pi\in\mathrm{S}_n\ .$$
The isomorphism $\mu$ translates this action to the action of $\mathrm{S}_n$ on the ring $\text{S}_\Z^{-1}\Diff (n)$. It turns out that the subring $\Diff(n)$ is preserved by this 
action. We present the formulas for the action of the generators $s_i$ of $\mathrm{S}_n$.
\begin{equation}\label{acofsn1}\begin{array}{l}{\displaystyle 
s_i(\Z^i) =- \Z^{i+1} \th_{i,i+1} \ ,\ s_i(\Z^{i+1}) = \Z^i \frac{1}{\th_{i,i+1}}\ ,\ s_i(\Z^j) =\Z^{j}\ \text{for}\ j\neq i,i+1\ ,}\\[.5em]
{\displaystyle s_i(\der_i) = - \frac{1}{\th_{i,i+1}} \der_{i+1}\ ,\ s_i(\der_{i+1}) = \th_{i,i+1} \der_i\ ,\ s_i(\der_j) =\der_{j}\ \text{for}\ j\neq i,i+1\ ,}\\[1.2em]
s_i(\th_j)=\th_{s_i(j)}\ .\end{array}\end{equation}

\paragraph{3.} For the $\R$-matrix description of the diagonal reduction algebra $\DRn$ in \cite{KO5} we used the ring $\Diff(n,N)$ formed by $N$ copies of the ring $\Diff(n)$. We do not know an analogue of the isomorphism $\mu$ for the ring $\Diff(n,N)$. However a straightforward analogue 
of the formulas (\ref{acofsn1}) provides an action of $\mathrm{S}_n$ by automorphisms on the ring $\Diff(n,N)$.

\vskip .2cm
We recall that the ring $\Diff(n,N)$ is a $\Uh$-bimodule with the generators $\Z^{j,\alpha}$ and $\der_{j,\alpha}$, $j=1,\dots,n$, $\alpha=1,\dots,N$.
The ring $\Diff (n,N)$ is free as a one-sided $\Uh$-module; the left and right $\Uh$-module structures are related by
\begin{equation}\label{relsweightsN}\th_i\Z^{j,\alpha}=\Z^{j,\alpha}(\th_i+\delta_i^j)\ ,\ \th_i\der_{j,\alpha}=\der_{j,\alpha}(\th_i-\delta_i^j)\ . \end{equation}
The defining relations for the generators $\Z^{j,\alpha}$ and $\der_{j,\alpha}$, $j=1,\dots,n$, $\alpha=1,\dots,N$, read
\begin{equation}\label{defhdiffAN}
 \Z^{i,\alpha} \Z^{j,\beta}=\sum_{k,l}\RR_{kl}^{ij}\Z^{k,\beta} \Z^{l,\alpha}\ ,\ \der_{i,\alpha}\der_{j,\beta}=\sum_{k,l}\RR_{ji}^{lk}\der_{k,\beta}\der_{l,\alpha}\ ,\ 
 \Z^{i,\alpha}\der_{j,\beta}=\sum_{k,l}\RR_{lj}^{ki}[\varepsilon_k]\der_{k,\beta} \Z^{l,\alpha}-\delta^{\alpha}_{\beta}\delta_{j}^i\ ,
\end{equation}
or, in components, 
\begin{equation}\label{rela1}
\Z^{i,\alpha} \Z^{j,\beta} = \frac{1}{\th_{ij}} \Z^{i,\beta} \Z^{j,\alpha} + \frac{\th_{ij}^2-1}{\th_{ij}^2} \Z^{j,\beta} \Z^{i,\alpha}\ ,\ 
\Z^{j,\alpha} \Z^{i,\beta} = - \frac{1}{\th_{ij}} \Z^{j,\beta} \Z^{i,\alpha} + \Z^{i,\beta} \Z^{j,\alpha}\ ,\ 1\leq i<j\leq n\ ,
\end{equation}
\begin{equation}\label{rela2}
\der_{i,\alpha} \der_{j,\beta} = -\frac{1}{\th_{ij}} \der_{i,\beta} \der_{j,\alpha} + \frac{\th_{ij}^2-1}{\th_{ij}^2} \der_{j,\beta} \der_{i,\alpha}\ ,\ 
\der_{j,\alpha} \der_{i,\beta} = \frac{1}{\th_{ij}} \der_{j,\beta} \der_{i,\alpha} + \der_{i,\beta} \der_{j,\alpha}\ ,\ 1\leq i<j\leq n\ ,
\end{equation}
\begin{equation}\label{rela3}
\Z^{i,\alpha} \der_{j,\beta} = \der_{j,\beta} \Z^{i,\alpha}\ ,\ \Z^{j,\alpha} \der_{i,\beta} = \frac{\th_{ij}(\th_{ij}+2)}{(\th_{ij}+1)^2} \der_{i,\beta} \Z^{j,\alpha}
\ ,\ 1\leq i<j\leq n\ ,
\end{equation}
\begin{equation}\label{rela4}
\Z^{i,\alpha} \der_{i,\beta} = \sum_{k=1}^n \frac{1}{1-\th_{ik}} \der_{k,\beta} \Z^{k,\alpha} - \delta^\alpha_\beta\ ,\ 1\leq i\leq n\\ .
\end{equation}
\begin{lemm}\hspace{-.2cm}. {\rm The maps $s_i$, $i=1,\dots,n-1$, defined on the generators of $\Diff(n,N)$ by 
\begin{equation}\label{acsnN1}
\begin{array}{l}{\displaystyle 
s_i(\Z^{i,\alpha}) = - \Z^{i+1,\alpha} \th_{i,i+1}\ ,\ s_i(\Z^{i+1,\alpha}) = \Z^{i,\alpha} \frac{1}{\th_{i,i+1}}\ ,\ s_i(\Z^{j,\alpha}) = \Z^{j,\alpha}\ 
\text{for}\ j\neq i,i+1\ ,}\\[.5em]
{\displaystyle 
s_i(\der_{i,\alpha}) = - \frac{1}{\th_{i,i+1}} \der_{i+1,\alpha}\ ,\ s_i(\der_{i+1,\alpha}) = \th_{i,i+1} \der_{i,\alpha}\ ,\ 
s_i(\der_{j,\alpha}) =\der_{j,\alpha}\ \text{for}\ j\neq i,i+1\ ,}\\[1.2em]
s_i(\th_j)=\th_{s_i(j)}\ ,\end{array}\end{equation}
extend to automorphisms of the ring $\Diff(n,N)$. Moreover, these automorphisms satisfy the Artin relations and therefore give the action of the  
symmetric group $\mathrm{S}_n$ by automorphisms.
}\end{lemm}
\noindent{\it{Proof.}} After the formulas (\ref{acsnN1}) are written down, the verification is a direct calculation. 
\hfill$\square$

\paragraph{4.} The diagonal reduction algebra $\DRn$ is a $\Uh$-bimodule with the generators $\teL_i^j$, $i,j=1,\dots,n$. The defining relations 
of $\DRn$ are given by the reflection equation, see \cite{KO5}
$$\RR_{12}\teL_1\RR_{12}\teL_1-\teL_1\RR_{12}\teL_1\RR_{12}=\RR_{12}\teL_1-\teL_1\RR_{12}\ ,$$
where $\teL=\{\teL_i^j\}_{i,j=1}^n$ is the matrix of generators (we refer to \cite{C,S,RS,KS,IO,IOP,IMO1,IMO2} for various aspects and applications 
of the reflection equation). 

\vskip .2cm
For each $N$ there is a homomorphism (\cite{KO5}, Section 4.1)
$$\tau_N\colon \DRn\to\Diff(n,N)\ \ \text{defined by}\ \ 
\tau_N(\teL_i^j)=\sum_\alpha \Z^{j,\alpha}\der_{i,\alpha}\ .$$
Moreover $\tau_N$ is an embedding for $N\geq n$. 

\vskip .2cm
The formulas (\ref{acsnN1}) show that the image of $\tau_N$ is preserved by the automorphisms $s_i$. The element $s_i(\tau_N(\teL^j_k))$ 
can be written by the same formula for all $N$. Since $\tau_N$ is injective for $N\geq n$ we conclude that the formulas (\ref{acsnN1})
induce the action of the symmetric group $\mathrm{S}_n$ on the diagonal reduction algebra $\DRn$ by automorphisms. 

\paragraph{5.} The isomorphism $\mu$ can be also used to translate the action (\ref{automq}) of the braid group by Zhelobenko operators to the 
action of the braid group by automorphisms on the ring $\mathfrak{k}[a_1,\dots,a_n]\otimes \text{T}^{-1}\text{W}_n$. It turns out that this action preserves 
the subring $\text{T}^{-1}\text{W}_n$. Moreover, let $\text{T}_0$ be the multiplicative set of $\text{T}$ generated by $X^jD_j-X^kD_k+\ell$, $1\leq j<k\leq n$, $\ell\in\mathbb{Z}$. Then the action of the operators $\q_i$, $i=1,\dots,n-1$, preserves the subring $\text{T}_0^{-1}\text{W}_n$. We present the formulas for the action of the operators $\q_i$, $i=1,\dots,n-1$:
$$\begin{array}{l}{\displaystyle 
\q_i(X^i) = \frac{1}{\mathcal{H}_{i,i+1}} X^{i+1}\ ,\ \q_i(X^{i+1}) = X^i \mathcal{H}_{i,i+1}\ ,\ \q_i(X^j)=X^j\ \text{for}\ j\neq i,i+1\ ,}\\[.5em] 
{\displaystyle \q_i(D_i) = D_{i+1} \mathcal{H}_{i,i+1}\ ,\ \q_i(D_{i+1}) = \frac{1}{\mathcal{H}_{i,i+1}} D_i\ ,\ \q_i(D_j)=D_j\ \text{for}\ j\neq i,i+1\ .}
\end{array}$$

\paragraph{6.} The isomorphism (\ref{fromdifftow}) allows  to construct a $2n$-parametric family of $\Diff(n)$-modules different from the highest weight representations. 
Let $\vec{\gamma}:=\{\gamma_1,\dots,\gamma_n\}$ 
be a sequence of length $n$ of complex numbers such that $\gamma_i-\gamma_j\notin\mathbb{Z}$ for all $i,j=1,\dots,n$, $i\neq j$. 
Let $V_{\vec{\gamma}}$ be the vector space with the basis 
$$v_{\vec{j}}:=(X^1)^{j_1+\gamma_1}(X^2)^{j_2+\gamma_2}\dots (X^n)^{j_n+\gamma_n}\ ,\ 
\text{where}\ \vec{j}:=\{ j_1,\dots,j_n\}\ ,\ j_1,\dots,j_n\in\mathbb{Z}\ .$$ 
Under the conditions on $\vec{\gamma}$, $V_{\vec{\gamma}}$ is naturally a $\text{T}^{-1}\text{W}_n$-module. Define the action of the elements $a_k$ on the space $V_{\vec{\gamma}}$ by $a_k\colon v_{\vec{j}}\mapsto A_k v_{\vec{j}}$ where 
$\vec{A}:=\{A_1,\dots,A_n\}$ is another sequence of length $n$ of complex numbers. Then $V_{\vec{\gamma}}$ becomes an 
$\mathfrak{k}[a_1,\dots,a_n]\otimes \text{T}^{-1}\text{W}_n$-module and therefore $\Diff(n)$-module which we denote by 
$V_{\vec{\gamma},\vec{A}}$. The central operator $c_k$ acts on $V_{\vec{\gamma},\vec{A}}$ 
by scalar multiplication on $A_k$.

\vskip .2cm
\noindent{\bf Acknowledgments.} The work of O. O. was supported by the Program of Competitive Growth of Kazan Federal University.

\end{document}